\documentclass[11pt]{article}
\usepackage[english]{babel}
\usepackage[hmarginratio=3:2]{geometry}
\usepackage[titletoc]{appendix}
\usepackage{authblk}
\usepackage{graphicx}
\usepackage{amsmath}
\usepackage{amsthm}
\usepackage{amssymb}
\usepackage[mathscr]{eucal}
\usepackage{textcomp}
\usepackage[latin1]{inputenc}
\usepackage{latexsym}
\usepackage{hyperref}
\usepackage[dvipsnames]{xcolor}
\usepackage{framed}
\usepackage{tikz-cd}
\usepackage{ifsym}
\usepackage{enumerate}
\usepackage{sectsty}
\usepackage{chngcntr}

\numberwithin{equation}{section}

\numberwithin{subsection}{section}

\allowdisplaybreaks[1]

\newtheorem{theorem}{Theorem}[section]
\newtheorem{lemma}[theorem]{Lemma}

\newcommand\cE{\mathcal{E}} \newcommand\cF{\mathcal{F}}
 
\newcommand\cI{\mathcal{I}} 
 \newcommand\cL{\mathcal{L}}
 \newcommand\cN{\mathcal{N}}

 \newcommand\cX{\mathcal{X}}
\newcommand\cY{\mathcal{Y}}

\newcommand\QQ{\mathbb{Q}}

 \newcommand\ZZ{\mathbb{Z}}

\newcommand\rK{\mathrm{K}}

\newcommand\arr{\ifinner\to\else\longrightarrow\fi}

\newcommand\arrto{\ifinner\mapsto\else\longmapsto\fi}

\def\displaytimes_#1{\mathrel{\mathop{\times}\limits_{#1}}}

\def\displayotimes_#1{\mathrel{\mathop{\bigotimes}\limits_{#1}}}

\newcommand\aut{\operatorname{Aut}}

\newcommand\rk{\operatorname{rk}}

\newlength{\ignora}

\newcommand{\mmu}{\boldsymbol{\mu}}

\newcommand{\GL}{\mathrm{GL}}

\renewcommand{\epsilon}{\varepsilon}

\renewcommand{\k}[1][*]{\operatorname{K}'_{#1}}

\newcommand{\R}{\mathrm{R}}

\newcommand{\kg}[2][*]{\rK'_{#1}(#2)_{\bf g}}

\newcommand{\kt}[2][*]{\rK'_{#1}(#2)_{\bf t}}

\newcommand{\ci}[1][\relax]{\cI_{\mmu_{{#1}}}}

\newcommand{\autm}[1][\infty]{\aut{\mmu_{#1}}}

\DeclareMathOperator{\Ind}{Ind}
\DeclareMathOperator{\Res}{Res}

\title{On a formula for the equivariant Euler characteristic of a $G$ sheaf}
\author[1]{Qiangru Kuang}
\affil[1]{\textit{Scuola Internazionale Superiore di Studi Avanzati,  Via Bonomea 265, 34136, Trieste, Italy}}
\author[2]{Francesco Sala}
\affil[2]{\textit{Scuola Normale Superiore, Piazza dei cavalieri 7, 56126, Pisa, Italy}}
\date{}

\begin{document}
\maketitle

\begin{abstract}
In the papers \cite{fischbacher2009equivariant}, \cite{kock}, H. Fischbacher-Weitz and B. K\"ock computed the equivariant Euler characteristic of a $G-$sheaf on a $G$-curve $X$ over a field. Using a form of the Riemann-Roch theorem for quotient stacks proved by the first author we extend their computations to the cases where $dim(X) >1$.
\end{abstract}

\bigskip
\bigskip

\section{Introduction}

In the papers \cite{fischbacher2009equivariant}, \cite{kock}, H. Fischbacher-Weitz and B. K\"ock computed the equivariant Euler characteristic of a $G-$sheaf on a $G$-curve $X$ over a field. Their results holds in the full generality where the order of $G$ can be divisible by the characteristic of the base field. Our aim is to generalize their result, and to archieve it we make use of a version of Riemann-Roch theorem for quotient stacks.\\

In \cite{Vistoli}, A. Vistoli proved a decomposition theorem for the equivariant K-theory of a $G$-variety $X$, where $G$ is a finite group. A priori, this result allows to calculate the equivariant Euler characteristic of any $G$-sheaf on $X$; in order to archieve that, however, we need to understand how his formula behaves under push-forward.\\

In \cite{To}, B. Toen provided a Riemann-Roch isomorphism from the K-theory of a quotient DM stack (with coefficient in $\QQ^{ab}$) to the étale K-theory of its inertia, and he proved that this map is covariant with respect to any push-forward. However his result only holds for the case where the stabilizers of the points have order prime with the characteristic of the base field.\\

In the recent \cite{S}, the second author proves a version of Toen's result which covers the case where the stabilizers of the points have order divisible by the characteristic of the base field. Using this result, we can compute explicitly the equivariant Euler characteristic of any $G$-sheaf on $X$, in Vistoli's context.\\

We briefly recall some of the main results of \cite{S}. First, we recall that, following \cite{dan-olsson-vistoli}, a stack is called tame if the stabilizer of any point is linearly reductive. 

In general, a morphism of stacks $f\colon\cX\rightarrow \cY$ is called relatively tame if the relative inertia $\cI_{\cY}\cX$ has linearly reductive geometric fibers (see \cite{dan-olsson-vistoli}). For such a morphism, the push-forward $f_*$ descends to K-theory. \\

Let $\cX$ be a quotient DM stack $[X/\GL_n]$ and $\cI_{\mu}\cX$ be its \emph{cyclotomic inertia}, defined as the union of the morphism stacks $\underline{Hom}(B\mmu_n,\cX)$; it comes with a finite projection $\rho\colon \cI_{\mu}\cX \arr \cX$. 

The components of $\cI_{\mu}\cX$ are indexed by the conjugacy classes of subgroups of $\GL_n$ which are isomorphic to some $\mmu_n$. We will call these subgroups \emph{dual cyclic}.\\

Moreover, following \cite{Schadeck}, we have a decomposition $K([X/\GL_n])\simeq \prod_\sigma K([X/\GL_n])_\sigma$ by localizations indexed again by classes of dual cyclic subgroups. This allows us to define the \emph{tautological part} of the K-theory of $\cI_{\mu}\cX$ as $\kt{\cI\cX}:=\prod_\sigma K(\kt{\cI_\sigma\cX})_\sigma$.\\

Finally, there is a \emph{twist map} $\alpha_\cX\colon \cI_{\mu}\cX\rightarrow \prod_\sigma \cI_{\sigma}\cX\times B\sigma$; in the presentation $\cI_\sigma\cX=[X^\sigma/C_{\GL_n}(\sigma)]$ the twist sends an equivariant sheaf $\cF$ on $X^\sigma$ to $\oplus \cF^i\otimes \chi^i$, where $\chi$ is a generator of the ring of characters $R(\sigma)$ and $\cF^i$ is the corresponding eigensheaf.\\

Then we have a map given by the composition
\[
\k(\cX)\xrightarrow{\rho_*^{-1}} \kt{\cI\cX}^{\autm}\xrightarrow{\alpha_{\cX,*}}\Bigl(\bigoplus_{\substack{r}}\kg{\ci[r]\cX}\otimes\QQ(\zeta_r) \Bigr)^{\autm}
\]
(the last arrow is the composition of $\alpha_{\cX,*}$ and the projections $\R\mu_r\arr\QQ(\zeta_r)$, corresponding to the observation that $\QQ(\zeta_r)\simeq R(\mu_r)_{\mu_r}$).\\

Let $r_*\colon \Bigl(\bigoplus_{\substack{r}}\kg{\ci[r]\cX}\otimes\QQ(\zeta_r) \Bigr)^{\autm}\arr \k(\cX) $ be the above composition. Then the map $\cL :=r_*^{-1}$ gives

\begin{theorem}\label{main}
Let $\cX$ be a quotient stack with finite cyclotomic inertia over a base scheme $A=Spec(R)$. Then the above map gives a
 \emph{Lefschetz-Riemann-Roch} isomorphism, which is covariant with respect to proper push-forwards of relatively tame morphisms:
\[
\cL\colon \k(\cX)\longrightarrow \Bigl(\bigoplus_{\substack{r}}\kg{\ci[r]\cX}\otimes\QQ(\zeta_r) \Bigr)^{\autm}
\]
\end{theorem}

If $\cX$ is regular we can express $\cL$ in a different way, much more feasible for calculations.

Let us consider the map $r^*$ given by the composite
\[
\k(\cX)\xrightarrow{\lambda_{-1}(\cN)^{-1}\cdot\rho^*} \kt{\cI\cX}^{\autm}\xrightarrow{\alpha_\cX^*}\Bigl(\bigoplus_{\substack{r}}\kg{\ci[r]\cX}\otimes\QQ(\zeta_r) \Bigr)^{\autm}
\]
where $\cN$ is the conormal sheaf for $\rho\colon\cI\cX\arr\cX$ and the last map is - obviously - composed with the projections $\R(\mmu_r)\arr\QQ(\zeta_r)$.\\

We have the following:

\begin{theorem}\label{comp}
Let $\cX$ be a regular stack with finite cyclotomic inertia over $A$. Suppose that $\cX$ is either
\begin{enumerate}
\item of the form $[X/G]$, where $X$ is a scheme and $G$ is a finite flat group scheme over $A$,
\item Deligne-Mumford
\item tame.

\end{enumerate}  

Then the composition $$r^*\circ r_*\colon \Bigl(\bigoplus_{\substack{r}}\kg{\ci[r]\cX}\otimes\QQ(\zeta_r) \Bigr)^{\autm}\arr \Bigl(\bigoplus_{\substack{r}}\kg{\ci[r]\cX}\otimes\QQ(\zeta_r) \Bigr)^{\autm}$$
is equal to the endomorphism that on each component $\kg{\ci[r]\cX}\otimes\QQ(\zeta_r)$ is a multiplication by the rational number $\frac{\phi(r)}{r}$.
\end{theorem}

In particular the Lefschetz-Riemann-Roch map, in the regular case, is equal to $\bigoplus_{\substack{d}}r^*\cdot \frac{d}{\phi(d)}$.\\

Let now take again \(X\) to be a smooth projective variety over a field equipped with an action of a finite group scheme \(G \subseteq Aut(X)\) of order \(n\), such that the cyclotomic inertia $\cI_\mu([X/G])$ is finite over $[X/G]$. Denote the stack quotient by \(\mathcal X = [X/G]\) and the scheme quotient \(Y\). The stack $\cX$ satisfies the hypothesis of the above Theorem, so given a \(G\)-equivariant vector bundle \(\mathcal E\) on \(X\)  we can use the above results to compute the class of  $\chi(\cE)$ in the representation ring of $G$.

Given a cyclic subgroup $\sigma$, we denote by \(X^\sigma\) the fixed loci by the subgroup, which inherits an action by the centraliser \(C(\sigma)\).

Then our result states that 
\begin{theorem}
In the above hypotheses (that is $\cI_\mu([X/G])$ is finite over $[X/G]$)
\[
  \chi_G(X, \mathcal E) = \bigoplus_\sigma \bigoplus_i \frac{\chi(A^{\sigma, i})}{\varphi(|\sigma|)} \frac{|\sigma|}{|C(\sigma)|}\cdot \Ind_\sigma^G \iota(x^i)
\]
where \(A^{\sigma, i}\) is a \(C(\sigma)\)-equivariant vector bundle on \(X^\sigma\) that can be given explicitly.
\end{theorem}

The result is obtained by computing the euler characteristics of \(\mathcal E\) in two ways via the Grothendieck-Riemann-Roch theorem for stacks:
\[
  \begin{tikzcd}
    K(\mathcal X) \ar[r, "\mathcal L_{\mathcal X}"] \ar[d, "\pi_*"] & \bigoplus_{r,\sigma \in \overline C_r}K_{\mathbf g}(\tilde I_\sigma \mathcal X)\otimes \QQ(\zeta_r) \ar[d] \\
    K(BG) \ar[r, "\mathcal L_{BG}"] & \bigoplus_{r,\sigma \in \overline C_r} K_{\mathbf g}(\tilde I_\sigma BG)\otimes \QQ(\zeta_r)
  \end{tikzcd}
\]
Since \(\mathcal X\) is smooth, the upper row can be identified with the composition
\[
  K(\mathcal X) \xrightarrow{\lambda_{-1}(\mathcal N^\vee)^{-1} \cdot \rho^*} K(I \mathcal X) \xrightarrow{m^*} K(\tilde I \mathcal X) \xrightarrow{\frac{r}{\phi(r)}} K(\tilde I \mathcal X),
\]
which by constuction lands in the tautological part of the \(K\)-theory of the inertia stack. On the other hand, \(\mathcal L_{BG}^{-1}\) is the composition
\[
  \bigoplus_{\sigma} K(BC(\sigma))_{\mathbf g} \otimes \tilde R\sigma \to \bigoplus_{\sigma} K(BC(\sigma)) \otimes R\sigma \xrightarrow{m_*} \bigoplus_{\sigma} K(BC(\sigma)) \xrightarrow{\Ind_{C(\sigma)}^G} K(BG).
\]

\section{Computation}

\subsection{Lower composition}

\paragraph{Inclusion of geometric part of \(BG\)}

The inclusion \(K_{\mathbf g}(BH) \to K(BH)\) sends the unit to \(\frac{kH}{H}\).

\paragraph{Inclusion \(\tilde R\sigma \to R\sigma\)}

Suppose \(\sigma\) is dual cyclic of order \(r\). We denote by \(\iota: \QQ[x]/\Phi_r(x) \to \QQ[x]/x^r - 1\) the section of the projection \(\QQ[x]/x^r - 1 \to \QQ[x]/\Phi_r(x)\).

\paragraph{Inertia stack of \(BG\) and representations}

Let \(\sigma\) be a dual cyclic subgroup of \(G\) of order \(r\) and let \(H = C(\sigma)\) be its centraliser. The twisting operation on the \(r\)th cyclotomic inertia stack is given by the pullback induced by multiplication \(m: H \times \sigma \to H\). The anti-twist is given by \(m_*\), whose effect is taking invariants relative to the subgroup \(\ker(m) = \{(x^{-1}, x): x \in \sigma\}\). In formula, for an \(H\)-representation \(V\) and a character \(\chi\) of \(\sigma\)
\begin{align*}
  m_*: RH \otimes R\sigma &\to RH \\
  V\otimes \chi &\mapsto \chi\text{-isotypical part of } \Res_\sigma^H V,
\end{align*}
which inherits a natural \(H\)-action.

Since \(H\) is the centraliser of \(\sigma\), we have $\Res_\sigma^H \Ind_\sigma^H \chi=\frac{|H|}{|\sigma|}\chi$ for any $\sigma$-character $\chi$ and
\[
  \Res_\sigma^H k[H] = \Res_\sigma^H\Ind_\sigma^H k[\sigma] = k[\sigma]^{\oplus |H/\sigma|}=\frac{|H|}{|\sigma|}\bigoplus_{\chi \in \widehat\sigma}\chi =\bigoplus_{\chi \in \widehat\sigma}\Res_\sigma^H \Ind_\sigma^H \chi.
\]
As a consequence we note that the $\chi-$isotypic part of $\Res_\sigma^H k[H]$ is exactly $\Ind_\sigma^H \chi$, whence
\[
  m_*(\frac{k[H]}{|H|} \otimes \chi) = \frac{1}{|H|} \Ind_\sigma^H \chi.
\]
In particular, the composition 
\[
\tilde R \sigma \simeq K(BH)_{\mathbf g} \otimes \tilde R \sigma \hookrightarrow K(BC(\sigma)) \otimes R\sigma \xrightarrow{m_*} K(BC(\sigma)) 
\]
is the same as
\[
\tilde R \sigma \simeq K(B\sigma)_{\mathbf g} \otimes \tilde R \sigma \hookrightarrow K(B\sigma) \otimes R\sigma \xrightarrow{m_*} K(B\sigma)\xrightarrow{\Ind_\sigma^{C(\sigma)}}  K(BC(\sigma)).
\]
up to a correcting factor of $\frac{r}{|H|}=\frac{|\sigma|}{|C(\sigma)|}$.\\

Finally, by \ref{comp}, the map 
\[
\tilde R \sigma \simeq K(B\sigma)_{\mathbf g} \otimes \tilde R \sigma \hookrightarrow K(B\sigma) \otimes R\sigma \xrightarrow{m_*} K(B\sigma)
\]
is equal to $\frac{1}{r}\cdot \iota\colon\mathbb{Q}(\zeta_r)\rightarrow \QQ[x]/x^r - 1 $.\\

Now we are ready to compute the lower composition. We index the inertia stack of \(BG\) by \(\overline{C} = \coprod \overline{C}_r\), where \(\overline{C}_r\) is the conjugacy classes of monomorphisms \(\mu_r \to G\):
\[
  IBG = \coprod_{\sigma \in \overline C_r} BC(\sigma).
\]
Then the morphism \(\mathcal L_{BG}\) is the composition
\[
  K(BG) \to \bigoplus_\sigma K(BC(\sigma)) \to \bigoplus_\sigma K(BC(\sigma)) \otimes R\sigma \to \bigoplus_\sigma K_{\mathbf g}(BC(\sigma)) \otimes \tilde R\sigma
\]
whose inverse on the summand labelled by \(\sigma \in \overline C_r\) given by
\[
  \begin{tikzcd}[row sep=small]
    K_{\mathbf g} (BC(\sigma))\otimes \tilde R \sigma \ar[r, hook] \ar[d, equal] & K(BC(\sigma)) \otimes R\sigma \ar[r, "m_*"] \ar[d, equal] & K(BC(\sigma)) \ar[r, "\Ind_{C(\sigma)}^G"] & K(BG)\\
    \QQ \otimes \QQ[x]/(\Phi_r(x)) & K(BC(\sigma)) \otimes \QQ[x]/x^r - 1 \\
    1 \otimes x \ar[r, mapsto] & \frac{kC(\sigma)}{|C(\sigma)|} \otimes \iota(x) \ar[r, mapsto] & \frac{r}{|C(\sigma)|}\cdot\frac{1}{r} \Ind_\sigma^{C(\sigma)} \iota(x) \ar[r, mapsto] & \frac{1}{|C(\sigma)|} \Ind_\sigma^G i(x)
  \end{tikzcd}
\]

\subsection{Upper composition}

The upper composition is considerably easier, thanks to \ref{comp} and the fact that the twisting map \(m^*\) is easier to describe. Write
\[
  I\mathcal X = \bigoplus_{r,\sigma \in \overline C_r}[X^\sigma/C(\sigma)].
\]
Given a \(G\)-equivariant vector bundle \(\mathcal E\) on \(X\), let \((\mathcal E_\sigma)_\sigma = \rho^*\mathcal E\) be its pullback along \(\rho: I \mathcal X \to \mathcal X\) and let \((N_\sigma)_\sigma\) be the normal bundle of \(\rho\). On each component of the inertia stack we split the virtual bundle
\[
  \frac{\mathcal E_\sigma}{\lambda_{-1}(N_\sigma^\vee)} := A^\sigma = \bigoplus_{i \in \hat \sigma} A^{\sigma, i}
\]
into isotypical components as \(\sigma\)-representations, obtained by restriction from \(C(\sigma)\). Here \(\hat \sigma\) is the character group of \(\sigma\), which can be identified with \(\ZZ/r\ZZ\).\\

The upper composition thus sends 
\[
\mathcal E \longrightarrow \bigoplus_{r,\sigma\in \overline C_r, i\in\hat\sigma} \frac{r\cdot \zeta_r^i}{\phi(r)}\cdot A_{\mathbf g}^{\sigma, i}
\]
where $A^{\sigma, i}_{\mathbf g}$  denotes the  projection to the geometric part.

\subsection{Total composition}

We have the vertical map
\[
\mathcal I\pi_*\colon \bigoplus_{r,\sigma\in \overline C_r} K_{\mathbf g}(I_\sigma\mathcal X)\otimes\QQ(\zeta_r)\longrightarrow \bigoplus_{r,\sigma\in \overline C_r} K_{\mathbf g}(I_\sigma BG)\otimes\QQ(\zeta_r)\simeq \bigoplus_{r,\sigma\in \overline C_r} \QQ(\zeta_r)
\]
which, since $\pi$ is representable, is the identity on the $\QQ(\zeta_r)$ components and sends $A_{\mathbf g}^{\sigma, i}$ to $\chi(A^{\sigma, i})$.

Indeed we have a commuting diagram
\[
\begin{tikzcd}
K(I_\sigma\mathcal X) \ar[r]\ar[d,"\pi_*"] & K_{\mathbf g}(I_\sigma\mathcal X) \ar[d , "I\pi_*"] \\
K(I_\sigma BG) \ar[r] & K_{\mathbf g}(I_\sigma BG) 
\end{tikzcd}
\]
the horizontal maps being projections to the geometric part; the bottom-left composition is then easily seen to be exactly the Euler characteristic map.\\

Combining this with the previous calculation of the lower composition, we immediately get:
\[
  \chi_G(X, \mathcal E) = \bigoplus_{r,\sigma \in \overline C_r} \bigoplus_{i \in \hat \sigma} \frac{\chi(A^{\sigma, i})}{\phi(r)} \frac{r}{|C(\sigma)|}\cdot\Ind_\sigma^G \iota(\zeta_r^i).
\]

\section{Explicit formula for curves}

In this section we specialise to the case where \(X\) is a smooth projective curve and $G$ is discrete. Then \(Y\) is also a smooth projective curve. In this case we aim to recover the results of \cite{fischbacher2009equivariant}, \cite{kock}. 

Let \(n = |G|\). First note that we can give more concrete description of the \(\chi(A^{\sigma, i})\)'s. Over the identity sector corresponding to \(\sigma = 1\) it is the holomorphic Euler characteristic of \(\mathcal E\). For \(\sigma \ne 1\), as the coarse moduli space of \([X^\sigma/C(\sigma)]\) is zero-dimensional, \(\chi(A^{\sigma, i})\) is nothing but the dimension of the virtual representation \(A^{\sigma, i}\).

Thus
\[
  \chi_G(X, \mathcal E) = \chi(X, \mathcal E) \frac{kG}{n} + \bigoplus_{r > 1,\sigma \in \overline C_r} \bigoplus_{i \in \hat \sigma} \frac{\dim A^{\sigma, i}}{\phi(r)}\frac{r}{|C(\sigma)|}\cdot \Ind_\sigma^G \iota(\zeta_r^i).
\]


Let us reindex the summation. We fix a representative for each \(\tilde \sigma \in \overline{C}_r\), which we also call \(\tilde \sigma\). Choose also a set of representatives \(\{\tilde x\}\) for each $G$-orbit with non-trivial stabilisers. The action being effective, for each $1\ne\tilde\sigma$ the fixed locus $X^\sigma$ is zero-dimensional, so we can regroup uniquely the $\{\tilde\sigma\}$'s to the sets $\{\tilde \sigma \subseteq G_{\tilde x}\}_{\tilde{x}}$. For each $x\in X$, denote by $e_x$ (resp. $e^t_x$) the ramification index (resp. the \emph{tame} ramification index). We have:
\begin{align*}
  \chi_G(X, \mathcal E)
  &= \chi(X, \mathcal E) \frac{kG}{n} + \sum_{\tilde x} \bigoplus_{1 \ne \tilde\sigma \subseteq G_{\tilde x}} \bigoplus_{i \in \hat \sigma} \frac{|X^\sigma|\cdot\dim A^{\sigma, i}}{\phi(r)} \frac{r}{|C(\sigma)|}\cdot\Ind_\sigma^G \iota(\zeta_r^i) \\
  &= \chi(X, \mathcal E) \frac{kG}{n} + \sum_{x \in X} \frac{e_x}{n} \Ind_{G_x}^G \bigoplus_{1 \ne \tilde\sigma \subseteq G_{x}} \bigoplus_{i \in \hat \sigma} \frac{|X^\sigma|\cdot\dim A^{\sigma, i}}{\phi(r)} \frac{r}{|C(\sigma)|}\cdot\Ind_\sigma^{G_x} \iota(\zeta_r^i).
\end{align*}

At this point we observe that making use of \ref{comp}, we can rewrite the expression
\[
\bigoplus_{1 \ne \tilde\sigma \subseteq G_{x}} \bigoplus_{i \in \hat \sigma} \frac{\dim A^{\sigma, i}}{\phi(r)} \frac{r}{|C(\sigma)|}\cdot\Ind_\sigma^{G_x} \iota(\zeta_r^i) \qquad (*)
\]
Indeed, let $A^\sigma= \Res^{G_x}_{C(\sigma)} V^\sigma $, where $V^\sigma$ is the $G_x$-virtual sheaf $\frac{\mathcal E_\sigma}{\lambda_{-1}(N_\sigma^\vee)}$.

Then we know from \ref{comp}, that the following composition of maps
\[
\begin{tikzcd}
R(G_x)\ar[r,"\Res"] & \underset{\tilde\sigma \subseteq G_{x}}{\underset{r,\sigma }{\bigoplus} }R(C_{G_x}(\sigma))\ar[r,"m^*"] &\underset{\tilde\sigma \subseteq G_{x}}{\underset{r,\sigma }{\bigoplus} } R(C_{G_x}(\sigma))_{\mathbf g}\otimes \QQ(\zeta_r)\\
\underset{\tilde\sigma \subseteq G_{x}}{\underset{r,\sigma }{\bigoplus} } R(C_{G_x}(\sigma))_{\mathbf g}\otimes \QQ(\zeta_r)\ar[r,"m_*"]& \underset{\tilde\sigma \subseteq G_{x}}{\underset{r,\sigma }{\bigoplus} } R(C_{G_x}(\sigma))\ar[r,"\Ind"] & R(G_x)
\end{tikzcd}
\]
given explicitly by
\[
\begin{tikzcd}
V^\sigma \ar[r] & \underset{\tilde\sigma \subseteq G_{x}}{\underset{r,\sigma }{\sum} } \Res^{G_x}_{C_{G_x}(\sigma)}(V^\sigma)\ar[r] & \underset{i\in\hat\sigma }{\underset{r,\sigma}{\sum}} \dim(V^{\sigma,i})\cdot\zeta^i_r\cdot\frac{r}{\phi(r)}\\
\underset{i\in\hat\sigma }{\underset{r,\sigma}{\sum}} \dim(V^{\sigma,i})\cdot\zeta^i_r\ar[r]\cdot\frac{r}{\phi(r)} & \underset{i\in\hat\sigma }{\underset{r,\sigma}{\sum}} \dim(V^{\sigma,i})\cdot \frac{\iota(\zeta^i_r)}{|C_{G_x}(\sigma)|}\cdot\frac{r}{\phi(r)}\ar[r]&
\underset{i\in\hat\sigma }{\underset{r,\sigma}{\sum}} \frac{\dim(V^{\sigma,i})}{\phi(r)}\frac{r}{|C_{G_x}(\sigma)|}\cdot \Ind^{G_x}_\sigma(\iota(\zeta^i_r)) 
\end{tikzcd}
\]
is the identity. Note that we computed the map $m_*$ exactly as we did at the end of Section 2.1.\\

In particular we can write
\[
\bigoplus_{1 \ne \tilde\sigma \subseteq G_{x}} \bigoplus_{i \in \hat \sigma} \frac{|X^\sigma|\cdot\dim A^{\sigma, i}}{\phi(r)} \frac{r}{|C(\sigma)|}\cdot\Ind_\sigma^{G_x} \iota(\zeta_r^i)=A_x - \iota_{\sigma = 1} \dim A_x
\]
noting that the missing term corresponding to $\sigma=1$ corresponds exactly to the geometric part of $A_x$.\\

\textbf{Remark:} We must have $\frac{|X^\sigma|}{|C(\sigma)|}=\frac{1}{|C_{G_x}(\sigma)|}$, since the connected components of the inertia stack corresponding to $1\ne\tilde\sigma\subseteq G_x$ are equal to those of $\mathcal IBG_x$. \\

Now we invoke the following lemma:

\begin{lemma}
  \label{lem:invert root of unity}
  Let \(H\) be a cyclic group of order \(m\) and \(\chi\) a non-trivial character. Then
  \[
    \frac{1}{1 - \chi} = - \frac{1}{m} \sum_{d = 1}^{m - 1} \chi^{-d}.
  \]
\end{lemma}

\begin{proof}
 Same as \cite{kock} Lemma 1.2.
\end{proof}

Noting that the action of $G_x$ on $ N_x^\vee$ factors through the tame part, we use the Lemma to write
\[
  \dim A_x = \dim \frac{\mathcal E_x}{1 - N_x^\vee} = \frac{1}{e^t_x} \dim \sum_{d = 1}^{e^t_x - 1} (-d) \mathcal E_x \otimes N_x^{-d} = -\frac{1}{e^t_x} \cdot \rk \mathcal E_x \frac{e^t_x(e^t_x - 1)}{2} = -\rk \mathcal E \frac{e^t_x - 1}{2}.
\]
and
\[
 A_x - \iota_{\sigma = 1} \dim A_x
  = A_x - \dim A_x \frac{kG_x}{e_x}
  = A_x + \rk \mathcal E \frac{e^t_x - 1}{2} \frac{kG_x}{e_x}.
\]
Putting all ingredients together,
\begin{align*}
  \chi_G(X, \mathcal E)
  &= \chi(X, \mathcal E) \frac{kG}{n} + \sum_{x \in X} \frac{e_x}{n} \Ind_{G_x}^G(A_x + \rk \mathcal E \cdot \frac{e^t_x - 1}{2} \frac{k G_x}{e_x}) \\
  &= \left(\chi(X, \mathcal E) + \frac{\rk \mathcal E}{2} \sum_x (e^t_x - 1) \right) \frac{kG}{n} + \sum_{x \in X} \frac{e_x}{n} \Ind_{G_x}^G \frac{\mathcal E_x}{1 - N_x^\vee}.
\end{align*}

Thus we have recovered \cite{fischbacher2009equivariant},
 Theorem 4.2. When the $G$-action is tame, we can recover the main theorem of  \cite{kock}.
Indeed one can use Hirzebruch-Riemann-Roch and Riemann-Hurwitz to express \(\chi(X, \mathcal E)\) in terms of rank and degree of \(\mathcal E\) and genus of \(Y\):
\[
  \chi(X, \mathcal E) = \deg \mathcal E + \rk \mathcal E (1 - g_X) = \deg \mathcal E + \rk \mathcal E \left(n (1 - g_Y) - \frac{1}{2} \sum_{x \in X} (e_x - 1)\right)
\]
and expand \(\frac{\mathcal E}{1 - N_x^\vee}\) using \ref{lem:invert root of unity}.

\end{document}